\documentclass[11pt]{article}
\usepackage{amssymb}
\usepackage{amsmath}

\newtheorem{theo}{Theorem}
\newtheorem{prop}{Proposition}
\newtheorem{lemm}{Lemma}
\newtheorem{coro}{Corollary}
\newtheorem{rema}{Remark}

\newcommand{\cqfd}
{%
\mbox{}%
\nolinebreak%
\hfill%
\rule{2mm}{2mm}%
\medbreak%
\par%
}
\newfont{\gothic}{eufb10}
\date{\empty}
\begin{document}
\title{On integral Hodge classes on uniruled or Calabi-Yau threefolds }
\author{Claire Voisin\\ Institut de math{\'e}matiques de Jussieu, CNRS,UMR
7586} \maketitle \setcounter{section}{-1}
\begin{flushright} {\it To Masaki Maruyama, on his 60th birthday}
\end{flushright}
\section{Introduction}

Let $X$ be a smooth complex projective variety of dimension $n$. The
Hodge conjecture is then true for rational Hodge classes of degree
$2n-2$, that is, degree $2n-2$ rational cohomology classes of Hodge
type $(n-1,n-1)$ are algebraic, which means that they are the
cohomology classes of algebraic cycles with
$\mathbb{Q}$-coefficients. Indeed, this follows from the hard
Lefschetz theorem, which provides an isomorphism:
$$\cup c_1(L)^{n-2}:H^2(X,\mathbb{Q})\cong H^{2n-2}(X,\mathbb{Q}),$$
from the fact that the isomorphism above sends the space of rational Hodge classes of degree $2$
onto the space of rational Hodge classes of degree $2n-2$, and from the Lefschetz theorem on
$(1,1)$-classes.

For integral Hodge classes, Koll\'ar \cite{kollar}, (see also \cite{voisou}) gave examples
 of  smooth complex projective manifolds
which do not satisfy the Hodge conjecture for integral degree $2n-2$ Hodge classes, for any $n\geq3$.
The examples are
smooth general  hypersurfaces $X$ of certain degrees in $\mathbb{P}^{n+1}$.
By the Lefschetz restriction theorem, such a variety satisfies
$$H^2(X,\mathbb{Z})=\mathbb{Z}H,\,H=c_1(\mathcal{O}_X(1)),$$
and
$$H^{2n-2}(X,\mathbb{Z})=\mathbb{Z}\alpha,\,<\alpha,H>=1.$$
Plane sections $C$ of $X$ have cohomology class $[C]=d\alpha,\,d=deg\,X$, because
$$deg\,C=d=<[C],H>.$$
Koll\'ar \cite{kollar} proves the following :
\begin{theo} \label{kollar}
Assume the degree $d$ of $X$ satisfies the property that $p^n$
divides $d$, for some integer $p$ coprime to $n!$. Then for a
general $X$, any curve $C$  in $X$ has degree divisible by $p$,
hence
 its cohomology class is a multiple of $p\alpha$. Thus
the class $\alpha$ is not algebraic, that is, it is not the
cohomology class of an algebraic cycle with integral coefficients.
\end{theo}
(Note that the condition that  $p$ be coprime to $n!$ can be
removed, using the degeneration used in section 4  of
\cite{voisou}.)

 The condition on the degree makes the canonical bundle of $X$ very ample, since
 the smallest possible degree available by this construction is $2^n$. It is thus natural
 to try to understand whether this is an artificial consequence of the method of construction, or
 whether  the positivity of the canonical bundle is essential.

 Another reason to ask whether one could find examples above
 with Kodaira dimension equal to $-\infty$ is the remark made in \cite{voisou} :

 \begin{lemm} \label{lemmeintro1}Let $X$ be a smooth rational complex projective manifold. Then
 the Hodge conjecture is true for integral Hodge classes of type $(n-1,n-1)$.
 \end{lemm}
(Note that the whole degree $2n-2$ cohomology of such an $X$ is of type
$(n-1,n-1)$, so the statement is that classes of curves generate $H^{2n-2}(X,\mathbb{Z})$
for a rational variety $X$.)

One can thus ask whether this criterion could be used to produce new examples of
unirational or rationally connected,  but non rational varieties
(we refer to \cite{CG}, \cite{B},\cite{I} for other criteria).
Namely, it would suffice to produce a smooth projective rationally connected
variety which does not satisfy the Hodge conjecture for
degree $2n-2$ integral cohomology classes.
The main result of this paper implies that in dimension $3$, this cannot be done:
\begin{theo}\label{main} Let $X$ be a smooth complex projective threefold which  either is
uniruled, or satisfies
$$K_X\cong \mathcal{O}_X,\,H^2(X,\mathcal{O}_X)=0.$$
 Then the Hodge conjecture is true for
integral degree $4$ Hodge classes on $X$.
\end{theo}
\begin{rema} {\rm Recall \cite{miyaoka} that a complex  projective threefold is uniruled,
that is swept out by rational curves, if and only if it has Kodaira dimension equal to $-\infty$.
Thus our condition is that either $\kappa(X)=-\infty$ or $K_X=\mathcal{O}_X$ and $H^2(X,\mathcal{O}_X)=0$.}
\end{rema}
Note that as an obvious corollary, we get the  following:
\begin{coro} Let $X$ be a smooth complex projective $n$-fold. Assume $X$ contains
a subvariety $Y$ which is
a smooth $3$-dimensional  complete intersection  of ample divisors, and satisfies one of the conditions in
Theorem \ref{main}. Then the Hodge conjecture is true for integral degree $2n-2$  Hodge classes
on $X$.
\end{coro}
Indeed, let $j$ be the inclusion of $Y$ into $X$. By Lefschetz restriction theorem, the map
$$j_*:H^4(Y,\mathbb{Z})=H_2(Y,\mathbb{Z})\rightarrow H_2(X,\mathbb{Z})=
H^{2n-2}(X,\mathbb{Z})$$ is an isomorphism. Thus the Hodge
conjecture for integral Hodge classes of degree $4$ on $Y$ implies
the Hodge conjecture for integral Hodge classes of degree $2n-2$ on
$X$. \cqfd Note that in higher dimensions, there are two possible
generalizations of the problem studied above. Namely, one can study
the Hodge conjecture for integral Hodge classes in degree $4$ or
$2n-2$. Both problems are birationally invariant, in the sense that
the two groups
$$Hdg^4(X,\mathbb{Z})/<[Z]>,\,Hdg^{2n-2}(X,\mathbb{Z})/<[Z]>,$$
where $<[Z]>$ denotes the subgroups generated by classes of
algebraic cycles with integral coefficients, are birational
invariants of a smooth complex projective manifold $X$ of dimension
$n$ (see \cite{voisou}). For both problems, it is clear that the
assumption ``uniruled'' will not be sufficient in higher dimension
to guarantee that the groups above vanish. Indeed, starting from one
of Koll\'ar's $3$-dimensional example $X,\alpha\in
H^4(X,\mathbb{Z})$, with $\alpha$ a non-algebraic integral Hodge class
(Theorem \ref{kollar}), we can consider the product
$$Y=X\times\mathbb{P}^1,$$
and both classes $$pr_1^*\alpha,\,pr_1^*\alpha\cup pr_2^*([pt])$$ in
degree $4$ and $6=2n-2$ respectively will give examples of
non-algebraic integral Hodge classes. However, one could ask whether
the analogue of Theorem \ref{main} holds for $X$ rationally
connected, and for integral Hodge classes of degree $4$ or $2n-2$ on
$X,\,n=dim\,X$.

The proof of  Theorem \ref{main} uses the Noether-Lefschetz locus for
surfaces $S$ in an adequately chosen ample linear system on $X$.
This  leads to simple criteria which guarantee that integral degree
$2$ cohomology classes on a given $S$ are generated over
$\mathbb{Z}$ by those which become algebraic on some small
deformation $S_t$ of $S$. The Lefschetz hyperplane section Theorem
allows then to conclude.

 In section \ref{sec1}, we state this criterion, which is an
 algebraic criterion concerning the infinitesimal variation of Hodge
 structure on $H^2(S)$, for varieties $X$ with
 $H^2(X,\mathcal{O}_X)=0$.
 In section \ref{sec2}, we prove that this  criterion is satisfied
for uniruled or $K$-trivial  varieties with trivial
$H^2(X,\mathcal{O}_X)$. In the case of $K$-trivial varieties, the
criterion had been also checked in \cite{voisinintjmath}, but the
proof given  here is substantially simpler. In section \ref{sec3}, a
refinement of this
 criterion for uniruled threefolds with $H^2(X,\mathcal{O}_X)\not=0$ is given
 and proven to hold for an adequate choice of linear system.

\vspace{0,5cm}

{\bf Thanks.} This work was started during the very interesting
conference ``{\it Arithmetic Geometry and Moduli Spaces}.
It is a pleasure to thank the organizers  for the nice atmosphere  they succeeded
to create.
I also wish to thank S. Mori for his help in the proof of Lemma \ref{KH2}
  and J. Starr for interesting discussions on related questions.

\section{An infinitesimal criterion \label{sec1}}
Let $X$ be a smooth complex projective $n$-fold. Let $j:S\hookrightarrow X$ be a surface
which is a smooth complete intersection of ample divisors. Thus by Lefschetz theorem, the
Gysin map:
$$j_*:H^2(S,\mathbb{Z})\rightarrow H^{2n-2}(X,\mathbb{Z})$$
is surjective.

We assume that the Hilbert scheme $\mathcal{H}$ of deformations of
$S$ in $X$ is smooth near $S$. This is the case  if $S$ is a smooth complete
intersection of sufficiently ample divisors. The space
$H^0(S,N_{S/X})$ is the tangent space to $\mathcal{H}$ at $S$. Let
$\rho:H^0(S,N_{S/X})\rightarrow H^1(S,T_S)$ be the Kodaira-Spencer
map, which  is the classifying map for the first order deformations
of the complex structure on $S$ induced
 by the universal family $\pi:\mathcal{S}\rightarrow\mathcal{H}$ of surfaces parameterized by
 $\mathcal{H}$.

 For $u\in H^1(S,T_S)$ we have the interior product with $u$:
 $$u\lrcorner: H^1(S,\Omega_S)\rightarrow H^2(S,\mathcal{O}_S).$$
 The criterion we shall use is the following:
 \begin{prop}\label{prop1} Assume there exists a $\lambda\in H^1(S,\Omega_S)$ such that
 the map
 \begin{eqnarray}\label{mulambda}\mu_\lambda:H^0(S,N_{S/X})\rightarrow H^2(S,\mathcal{O}_S),\\ \nonumber
 \mu_\lambda(n)=\rho(n)\lrcorner\lambda,
 \end{eqnarray}
 is surjective. Then any class $\alpha\in H^{2n-2}(X,\mathbb{Z})$ is algebraic.
 \end{prop}
 \begin{rema}{\rm Our assumptions imply immediately that the cohomology
 $H^{2n-2}(X,\mathbb{C})$ is of type $(n-1,n-1)$. Indeed, this last fact is equivalent to the vanishing
 of the space
 $H^{n}(X,\Omega_X^{n-2})$. On the other hand, interpreting the map $\mu_\lambda$ above
 in terms of infinitesimal variations of Hodge structures on the degree $2$ cohomology
 of the surfaces $S_t$ parameterized by $\mathcal{H}$, one sees that
 $Im\,\mu_\lambda$ is contained in
 $$Ker\,(j_*:H^2(S,\mathcal{O}_S)\rightarrow H^{n}(X,\Omega_X^{n-2})).$$
 Thus the assumptions imply that this last map $j_*$ is $0$, and as it is surjective by Lefschetz
 theorem,
 it follows that $H^{n}(X,\Omega_X^{n-2})=0$.
 }
 \end{rema}
\begin{rema}{\rm The assumption of Proposition \ref{prop1} is exactly the assumption
of Green's infinitesimal  criterion for the density of the Noether-Lefschetz locus
(see \cite{voisinbook}, 5.3.4),
which allows to conclude that real degree $2$ cohomology classes on $S$ can be
approximated by rational algebraic cohomology classes on nearby fibers $S_t$.
It had been already used in \cite{voisinens}, \cite{voisinduke} to construct interesting algebraic cycles
on Calabi-Yau threefolds.
 }
 \end{rema}
{\bf Proof of Proposition.} We refer to \cite{voisinbook}, chapter 5, for more details on infinitesimal variations of
Hodge structures. On a simply connected neighborhood $U$ in
$\mathcal{H}$ of the point $0\in\mathcal{H}$ parameterizing $S\subset X$, the restricted family
$$\pi:\mathcal{S}_U\rightarrow U$$
is differentiably trivial, and in particular the local system
$$H^2_\mathbb{Z}:=R^2\pi_*\mathbb{Z}/torsion$$
is trivial. Thus the locally free sheaf
$$\mathcal{H}^2:=H^2_\mathbb{Z}\otimes \mathcal{O}_U$$
is canonically trivial, and denoting by ${H}^2$ the corresponding
vector bundle on $U$, we get a canonical isomorphism
$${H}^2\cong U\times H^2(S,\mathbb{C}),$$
since the fiber of ${H}^2$ at $0$ is canonically isomorphic to $H^2(S,\mathbb{C})$.
Composing with the second projection gives us a
holomorphic map
$$\tau: H^2\rightarrow H^2(S,\mathbb{C}),$$
which on each fiber $H^2_t=H^2(S_t,\mathbb{C})$
is the natural identification
$$H^2(S_t,\mathbb{C})\cong H^2(S,\mathbb{C}).$$
Next the vector bundle $H^2$ contains a holomorphic subbundle
$F^1H^2$, which at the point $t\in U$ has for fibre the
subspace
$$F^1H^2(S_t):=H^{2,0}(S_t)\oplus H^{1,1}(S_t)\subset H^2(S_t,\mathbb{C}).$$
We shall denote by
$$\tau_1: F^1H^2\rightarrow  H^2(S,\mathbb{C})$$
the restriction of $\tau$ to $F^1H^2$.

The key point is the following fact, for which we refer to \cite{voisinbook}??:
\begin{lemm} For $\lambda\in H^{1}(S_t,\Omega_{S_t})$, choose any lifting
$\tilde{\lambda}\in F^1H^2_t$ of $\lambda$.
Then the surjectivity of the map
$$\mu_\lambda:H^0(S_t,N_{S_t/X})\rightarrow H^2(S_t,\mathcal{O}_{S_t})
$$
is equivalent to the fact that the map $\tau_1$ is a submersion at $\tilde{\lambda}$.
\end{lemm}
Having this, we conclude as follows: First of all, we observe that the assumption
of Proposition \ref{prop1} is
a Zariski open condition on $\lambda\in H^1(S,\Omega_S)$. Now, the space
$H^1(S,\Omega_S)=H^{1,1}(S)$ has a real structure, namely
$$H^{1,1}(S)=H^{1,1}(S)_{\mathbb{R}}\otimes\mathbb{C},$$
where $H^{1,1}(S)_{\mathbb{R}}=H^{1,1}(S)\cap H^2(S,\mathbb{R})$.
Thus if the assumption is satisfied for one $\lambda\in H^{1,1}(S)$,
it is satisfied for one real $\lambda\in H^{1,1}(S)_{\mathbb{R}}$.

In the Lemma above, choose for lifting $\tilde{\lambda}$ the class $\lambda$ itself.
Thus $\tilde{\lambda}$ is real, and so is $\tau_1(\tilde{\lambda})$.
As the assumption on $\lambda$ and the Lemma imply that $\tau_1$ is a submersion at
$\tilde{\lambda}$, so is the restriction
$$\tau_{1,\mathbb{R}}:H^{1,1}_{\mathbb{R}}\rightarrow H^2(S,\mathbb{R})$$
of $\tau_1$ to $\tau_1^{-1}(H^2(S,\mathbb{R}))$. Here we identified
$\tau_1^{-1}(H^2(S,\mathbb{R}))$ to
$$\cup_{t\in U}F^1H^2(S_t)\cap H^{2}(S_t,\mathbb{R})=
\cup_{t\in U}H^{1,1}(S_t)_{\mathbb{R}}=:H^{1,1}_{\mathbb{R}}.$$
As $\tau_{1,\mathbb{R}}$ is a submersion at $\tilde{\lambda}$, and $H^{1,1}(S)_{\mathbb{R}}$
is a smooth real manifold, because it is a real vector bundle on  $U$ and $U$ is smooth,
$Im\,\tau_{1,\mathbb{R}}$ contains an open set of
$H^2(S,\mathbb{R})$. On the other hand $Im\,\tau_{1,\mathbb{R}}$ is a cone.
We use now the following elementary Lemma:
\begin{lemm}\label{lemmecone} Let $V_\mathbb{Z}$ be a lattice, and
let $C$ be an open cone in $V_\mathbb{R}:=V_\mathbb{Z}\otimes\mathbb{R}$.
Then $V_\mathbb{Z}$ is generated over $\mathbb{Z}$ by the points in $V_\mathbb{Z}\cap C$.
\end{lemm}
We apply Lemma \ref{lemmecone} to
$V_\mathbb{Z}=H^2(S,\mathbb{Z})/torsion$ and to
$C=Im\,\tau_{1,\mathbb{R}}$. Thus we conclude that
$H^2(S,\mathbb{Z})/torsion$ is generated over $\mathbb{Z}$ by
classes $\alpha\in Im\,\tau_{1,\mathbb{R}}$. But by definition of
$\tau_1$, if $\alpha=\tau_{1,\mathbb{R}}(\lambda_t)$, for some
$\lambda_t\in H^{1,1}(S_t)_\mathbb{R}$, the corresponding class
$\alpha_t\in H^2(S_t,\mathbb{Z})/torsion$ identifies to $\lambda_t$
in $H^2(S_t,\mathbb{R})$. Thus the class
$$\lambda_t\in H^{1,1}(S_t)_\mathbb{R}\subset H^2(S_t,\mathbb{R})$$ is in fact
integral, hence it is  algebraic by Lefschetz theorem on
$(1,1)$-classes.

The conclusion is that, under the assumptions of Proposition \ref{prop1},
the lattice $H^2(S,\mathbb{Z})/torsion$ is generated over $\mathbb{Z}$ by integral classes which become
algebraic (i.e. are the class of a divisor) on some nearby fiber $S_t$.
As the torsion of $H^2(S,\mathbb{Z})$ is algebraic, the same conclusion holds for
$H^2(S,\mathbb{Z})$.

Finally, as the map $j_*:H^2(S,\mathbb{Z})\rightarrow H^{2n-2}(X,\mathbb{Z})$
is surjective, we conclude that $H^{2n-2}(X,\mathbb{Z})$ is generated
over $\mathbb{Z}$ by classes of
$1$-cycles in $X$.

\cqfd
\section{Proof of  Theorem \ref{main}\label{sec2} when $H^2(X,\mathcal{O}_X)=0$ }
In this section, we assume that $H^2(X,\mathcal{O}_X)=0$ and $X$
either has trivial canonical bundle or is uniruled.

In case where $X$ is uniruled, we have the following result:
\begin{lemm}\label{KH2} Let $X$ be a uniruled threefold. Then a smooth  birational
model $X'$ of $X$ carries an ample line bundle $H$ such that
$$H^2K_{X'}<0.$$
\end{lemm}
{\bf Proof.} As $X$ is uniruled, $X$ is birationally equivalent to
a $\mathbb{Q}$-Gorenstein threefold $Y$ which is either a singular Fano threefold,
or a Del Pezzo fibration over a smooth curve, or a conic bundle over
a $\mathbb{Q}$-Gorenstein surface.
Let us first prove the existence of an ample line bundle $H_Y$ on
$Y$ such that $K_YH_Y^2<0$:

a) If $Y$ is Fano, $-K_Y$ is ample, so we can take for $H_Y$ an integral multiple
of $-K_Y$.

b) Otherwise there is a morphism
$$\pi:Y\rightarrow B,$$
where $B$ is $\mathbb{Q}$-Gorenstein of dimension $1$ or $2$, and the relative canonical bundle
$K_\pi$ has the property that $-K_\pi$ is a relatively ample $\mathbb{Q}$-divisor.
Let $H_B$ be an ample line bundle on $B$, and choose for
$H_Y$ the $\mathbb{Q}$-divisor
$$H_Y=\pi^*H_B-\epsilon K_\pi,$$
where $\epsilon$ is a small rational number.
As $-K_\pi$ is relatively ample, $H_Y$ is ample for small enough $\epsilon$.
We compute now:
$$ H_Y^2K_Y=(\pi^*H_B-\epsilon K_\pi)^2(\pi^*K_B+K_\pi)$$
$$=\pi^*H_B^2K_\pi-2\epsilon K_\pi\pi^*H_B(\pi^*K_B+K_\pi)+O(\epsilon^2).$$
If $dim\, B=2$, the term $\pi^*H_B^2K_\pi$ is negative, so that for
small $\epsilon$, $ H_Y^2K_Y<0$. If $dim\, B=1$, the first term
vanishes but the second term is equal to $-2\epsilon
K_\pi^2\pi^*H_B$ and this is negative because $-K_\pi$ is relatively
ample.

Let now $Y, \,H_Y$ be as above, and let $\tau:X'\rightarrow Y$ be a desingularization of
$Y$. Thus $X'$ is a smooth birational model of
$X$. Then there is a relatively ample
divisor $E$ on $X'$ which  is supported on the exceptional divisor
of $\tau$.
Consider the $\mathbb{Q}$-divisor
$$H=\tau^*H_Y+\epsilon E,$$
for $\epsilon$ a sufficiently small rational number.
Then we have $K_{X'}=\tau^*K_Y+F$ where $F$ is supported on the exceptional divisor of $\tau$. This gives
$$ H^2K_{X'}=(\tau^*H_Y+\epsilon E)^2(\tau^*K_Y+F).$$
As $\tau^*H_Y^2F=0$, the dominating term is equal to
$\tau^*H_Y^2\tau^*K_Y=H_Y^2K_Y<0$. Thus for small $\epsilon$ we have
$H^2K_{X'}<0$.

\cqfd
From now on, we will, in the uniruled case, consider  $X'$ instead of $X$, which can be done since the
statement of  Theorem \ref{main} is invariant
under birational equivalence, and we will assume that $H$ satisfies the conclusion of Lemma \ref{KH2}.

Our aim in this section is to prove the following Proposition, which
by Proposition \ref{prop1} implies Theorem \ref{main} for uniruled
 and Calabi-Yau threefolds $X$ with $H^2(X,\mathcal{O}_X)=0$.
\begin{prop}\label{propmainsec2} Let $X$ be a smooth projective
uniruled or Calabi-Yau threefold such that $H^2(X,\mathcal{O}_X)=0$.
Let $H$ be an ample line bundle on $X$. In the uniruled case, assume
that $H$ satisfies $H^2K_X<0$. Then for $n$ large enough, and for
$S$ a generic surface in $\mid nH\mid$, there is a $\lambda\in
H^1(S,\Omega_S)$ which satisfies the property that
$$\mu_\lambda:H^0(S,\mathcal{O}_S(nH))\rightarrow
H^2(S,\mathcal{O}_S)$$ is surjective.
\end{prop}
To see that this is a reasonable statement, note that in the
$K$-trivial case, the spaces $ H^0(S,\mathcal{O}_S(nH))$ and $
H^2(S,\mathcal{O}_S)$ have the same dimension, since, if  $S\in\mid
n H\mid$, we have by adjunction
$$H^0(S,K_S)=H^0(S,\mathcal{O}_S(S))=H^0(S,N_{S/X}),$$
with $H^0(S,K_S)=H^2(S,\mathcal{O}_S)^*$. Thus the two spaces
involved in Proposition \ref{propmainsec2} have the same dimension.
In the uniruled case, we have:
\begin{lemm}  \label{12declem} Assume $X,H$ satisfies $H^2K_X<0$, then for
$S\in \mid nH\mid$, we have
$$h^0(\mathcal{O}_S(S))=h^0(K_S)+\phi(n),$$
where $\phi(n)=\alpha n^2+o(n^2),\,\alpha>0$.
\end{lemm}
{\bf Proof.} We have $K_S=K_X(S)_{\mid S}$. Thus
$$\chi(\mathcal{O}_S(S))=\chi(K_S(-K_X))$$
$$=\chi({K_X}_{\mid S})=\chi(\mathcal{O}_S)+\frac{1}{2}(K_{X\mid S}^2-K_{X\mid S}(K_{X\mid S}+nH_{\mid S}))$$
$$=\chi({K}_S)+\frac{1}{2}(nHK_X^2-nHK_X(nH+K_X)).$$
It follows that
$$\chi(\mathcal{O}_S(S))-\chi({K}_S)=-\frac{1}{2}n^2H^2K_X+{\rm affine\,\,linear\,\,\,term\,\,in}\,\,\,n.$$
On the other hand, for large $n$, the ranks
$$h^1(\mathcal{O}_S(S))=h^2(\mathcal{O}_X),\,
h^2(\mathcal{O}_S(S))=h^3(\mathcal{O}_X),$$
$$h^1(K_S)=h^1(\mathcal{O}_S)=h^1(\mathcal{O}_X),\,h^2(K_S)=\mathbb{C}$$
 do not depend on $n$.
It follows that we also have
$$h^0(\mathcal{O}_S(S))-h^0({K}_S)=-\frac{1}{2}n^2H^2K_X+{\rm affine\,\,linear\,\,\,term\,\,in}\,\,\,n,$$
which proves the result with $\alpha=-\frac{1}{2}H^2K_X>0$. \cqfd By
this Lemma, we conclude that in the $K$-trivial case and in the
uniruled case, we can assume that we have for $n$ large enough, and
$S\in\mid nH\mid$,
$$h^0(N_{S/X})=h^0(S,\mathcal{O}_S(S))\geq h^0(K_S)=h^2(\mathcal{O}_S).$$
This makes possible the surjectivity of the map
$$\mu_\lambda:H^0(S,N_{S/X})\rightarrow H^2(S,\mathcal{O}_S)$$
of (\ref{mulambda}), and also says that $\mu_\lambda$ is surjective
if and only if it has maximal rank.

Another way to see this is to introduce
$$V:=H^0(S,K_S),\,V':=H^0(S,N_{S/X}).$$
The bilinear map
\begin{eqnarray}\label{mu}\mu:V\times V'\rightarrow H^1(S,\Omega_S),\\ \nonumber
\mu(v,v')=\rho(v')\lrcorner v
\end{eqnarray}
and Serre's duality $H^1(S,\Omega_S)\cong H^1(S,\Omega_S)^*$ give
a dual map
$$q=\mu^*:H^1(S,\Omega_S)\rightarrow( V\otimes V')^*=H^0(\mathbb{P}(V)\times\mathbb{P}(V'),\mathcal{O}(1,1)),$$
given by
$$q(\lambda)(v\otimes v')=<\lambda,\mu(v\times v')>.$$
As we have
$$<\lambda,\rho(v')\lrcorner v>=-<\rho(v')\lrcorner\lambda,v>,$$
where the $<,>$ stand for Serre's duality on $H^1(S,\Omega_S)$ on the left and between
$H^0(S,K_S)$ and $H^2(S,\mathcal{O}_S)$ on the right, we see that
$q(\lambda)$ identifies to $\mu_\lambda\in Hom\,(V,V'^*)$.

Thus the condition that $\mu_\lambda$ has maximal rank for generic $\lambda$ is equivalent to the condition
that the hypersurface of
$\mathbb{P}(V)\times\mathbb{P}(V')$ defined by $q(\lambda)$ is non singular.

We shall use the following criterion:
\begin{lemm} \label{lemdim}
Given $\mu$ as in (\ref{mu}),  the generic hypersurface defined by
$q(\lambda)$ is non singular if the following set
\begin{eqnarray}
\label{Z} Z=\{(v_1,v_2)\in\mathbb{P}(V)\times\mathbb{P}(V'),\,\mu(v_1\times v_2)=0\,\,{\rm in}\,\,
H^1(S,\Omega_S)\}
\end{eqnarray}
satisfies
$$dim\,Z<dim\,\mathbb{P}(V').$$
\end{lemm}
{\bf Proof.} Assume to the contrary that the generic $q(\lambda)$ is singular.
Let
$$Z'\subset \mathbb{P}(H^1(S,\Omega_S))\times \mathbb{P}(V),$$
$$Z'=\{(\lambda,v),\,q(\lambda)\,\,{\rm is\,\,singular\,\,at\,\,} (v,v')\,\,{\rm for\,\,some}
\,\,v'\in \mathbb{P}(V')\}.$$
By assumption $Z'$ dominates $\mathbb{P}(H^1(S,\Omega_S))$. Clearly there is only one irreducible
component $Z'_g$ of $Z'$ which dominates $\mathbb{P}(H^1(S,\Omega_S))$.
Let $Z'_1$ be the second projection of
$Z'_g$ in $\mathbb{P}(V)$.

As $Z'_g$  dominates $\mathbb{P}(H^1(S,\Omega_S))$ we have
$$dim\,Z'_g\geq rk\,H^1(S,\Omega_S)-1.$$
On the other hand, the  fiber of
$Z'_g$ over the generic point $v_g$ of $Z'_1$ is equal
to
$$\mu(v_g\times V')^\perp.$$
Thus we have
$$dim\,Z'_g=dim\,Z'_1+rk\,H^1(S,\Omega_S)-1-rk\,\mu_{v_g},$$
where $\mu_{v_g}:V'\rightarrow H^1(S,\Omega_S)$ is the map
$v'\mapsto\mu(v_g\times v')$.

The condition $dim\,Z'_g\geq rk\,H^1(S,\Omega_S)-1$ is thus equivalent to
\begin{eqnarray}\label{ineg}dim\,Z'_1\geq rk\,\mu_{v_g}.
\end{eqnarray}
But on the other hand, the unique irreducible component
 $Z_0$ of $Z'_1$ which dominates $Z'_1$ has dimension
  equal to
$dim\,Z'_1+dim\,\mathbb{P}(V')-rk\,\mu_{v_g}$ and  inequality (\ref{ineg})
implies that this is $\geq dim\,\mathbb{P}(V')$.
\cqfd

Our first task will be thus to study the set $Z$ introduced in (\ref{Z}).
To this effect, we degenerate the surface $S\in \mid nH\mid$ to a surface with many nodes.
The reason for doing that is the following fact (cf \cite{voisinintjmath}):
\begin{lemm} \label{voisinlemme}Let $\mathcal{S}\rightarrow \Delta$ be a Lefschetz degeneration of surfaces
$S_t\in \mid nH\mid$, where the central fiber has ordinary double points
$x_1,\ldots,x_N$ as singularities.
Then the
limiting space
$$
\lim_{t\rightarrow0}
Im\,(q_t: H^1(S_t,\Omega_{S_t})\rightarrow ( H^0(S_t,K_{S_t})\otimes H^0(S_t,\mathcal{O}_{S_t}(nH)))^*)$$
which is a subspace of
$(H^0(S_0,K_{S_0})\otimes H^0(S_0,\mathcal{O}_{S_0}(nH)) )^*$
contains
for each $i=1,\ldots,N$ the multiplication-evaluation map which is the composite:
$$H^0(S_0,K_{S_0})\otimes H^0(S_0,\mathcal{O}_{S_0}(nH))
\rightarrow H^0(S_0,K_{S_0}(nH))\rightarrow H^0(K_{S_0}(nH)_{\mid x_i}).$$
\end{lemm}
To get surfaces with many nodes, we use  discriminant surfaces as in
\cite{barth}. We assume here that $H$ is very ample on $X$, and we
consider a generic symmetric $n$ by $n$ matrix $A$ whose entries
$A_{ij}$ are in $H^0(X,\mathcal{O}_X(H))$. Let
$\sigma_A:=discr\,A\in H^0(X,\mathcal{O}_X(nH))$ and $S_A$ be the
surface defined by $\sigma_A$.
\begin{theo} \label{barth}(Barth \cite{barth}) The surface $S_A$ has $N$ ordinary double points as singularities,
with
$$N=\begin{pmatrix}{n+1}\\{3}\end{pmatrix}H^3.$$
\end{theo}
Note that for large $n$, this grows like $\frac{n^3}{6}H^3$ while both dimensions
$h^0(K_S),\,h^0(\mathcal{O}_S(nH))$ of the spaces $V,\,V'$ grow like $h^0(\mathcal{O}_X(nH))$, that is like
$\frac{n^3}{6}H^3$ by Riemann-Roch.

Next we have the following Proposition, which might well be  known
already, but for which we could not find a reference:
\begin{prop} \label{propvan} Let $X$ be a smooth projective
threefold, and $H$ a very  ample line bundle on $X$ which satisfies
the property that $H^i(X,\mathcal{O}_X(lH))=0$ for $i>0,\,l>0$. Let
$S_A\in\mid nH\mid$ be a generic discriminant surface as above, and
let $W\subset X$ be its singular set. Then the cohomology group
$H^1(X,\mathcal{I}_W((n+2)H))$ vanishes.
\end{prop}
{\bf Proof.} Let $G=Grass(2,n)$ be the Grassmannian of
$2$-dimensional vector subspaces of $K:=\mathbb{C}^n$. The matrix
$A$ as above can be seen as a family of quadrics $A_x$ on
$\mathbb{P}(K)$ parameterized by $x\in X$, the surface $S_A$
corresponds to singular quadrics and the singular set $W$
parameterizes quadrics of rank $n-2$. Thus $W$ is via the second
projection in one-to-one correspondence with the following algebraic
set:
$$\widetilde{W}:=\{(l,x)\in G\times X, A_x\,\,{\rm is\,\,singular\,\,along}\,\,l\}.$$
 Let  $\mathcal{E}$ be the tautological rank $2$ quotient bundle
on $G$, whose fiber at $l$ is $H^0(\mathcal{O}_{\Delta_l}(1))$.
$\mathcal{E}$ is a quotient of $K^*\otimes\mathcal{O}_G$, and there
is the natural  map
$$e:S^2K^*\otimes\mathcal{O}_G\rightarrow K^*\otimes\mathcal{E}.$$
Let \begin{eqnarray} \label{F2dec}\mathcal{F}:=Im\,e.
\end{eqnarray} Clearly, a quadric $A \in S^2K^*$ on $\mathbb{P}(K)$
is singular along $\Delta_l$ if and only if it vanishes under the
map $e$ at the point $l$. Thus the set $\widetilde{W}$ is  the zero
locus of a section of the vector bundle
$$\mathcal{F}\boxtimes \mathcal{O}_X(H)$$
which is of rank
$2n-1$ on $G\times X$.
Note that the cokernel of $e$ identifies to
$\bigwedge^2\mathcal{E}=:\mathcal{L}$ where $\mathcal{L}$ is the Pl\"ucker line bundle on $G$.
Thus we have an exact sequence
\begin{eqnarray}
\label{exactF}0\rightarrow \mathcal{F}\rightarrow K^*\otimes\mathcal{E}\rightarrow \mathcal{L}\rightarrow0
\end{eqnarray}
on $G$.

As $\widetilde{W}$ is the zero set of a transverse section of a rank $2n-1$ vector bundle
on $G\times X$, its ideal sheaf admits the Koszul resolution:
$$0\rightarrow\bigwedge^{2n-1}\mathcal{F}^*\boxtimes \mathcal{O}_X((-2n+1)H)\rightarrow\ldots
\rightarrow \mathcal{F}^*\boxtimes \mathcal{O}_X(-H)\rightarrow
\mathcal{I}_{\widetilde{W}}\rightarrow 0.$$ Thus the space
$H^1(X,\mathcal{I}_W((n+2)H))=H^1(G\times
X,\mathcal{I}_{\widetilde{W}}\otimes pr_2^*((n+2)H))$ is the abutment of a
spectral sequence whose $E_1$-term is equal to
$$H^i(G\times X,\bigwedge^i\mathcal{F}^*\boxtimes\mathcal{O}_X((n+2-i)H)),\,i\geq1.$$

By K\"unneth decomposition and the vanishing assumptions, these
spaces split as:

$$H^i(G,\bigwedge^i\mathcal{F}^*)\otimes
H^0(X,(n+2-i)H),\,n+2-i>0,$$
$$H^i(G,\bigwedge^i\mathcal{F}^*)\otimes H^0(X,\mathcal{O}_X)\oplus
 H^{i-1}(G,\bigwedge^i\mathcal{F}^*)\otimes H^1(X,\mathcal{O}_X)$$
 $$\oplus H^{i-2}(G,\bigwedge^i\mathcal{F}^*)\otimes H^2(X,\mathcal{O}_X)
 \oplus H^{i-3}(G,\bigwedge^i\mathcal{F}^*)\otimes
 H^3(X,\mathcal{O}_X),\,i=n+2,$$
$$H^{i-3}(G,\bigwedge^i\mathcal{F}^*)\otimes H^3(X,(n+2-i)H),\,n+2-i<0.
$$
The proof of the proposition is thus concluded by the following
Lemma, which implies that the $E_1$-terms of the spectral sequence
above all vanish.
\begin{lemm}\label{vanlemm3dec} On the Grassmannian $G=Grass(2,n)$,
the bundle $\mathcal{F}$ being defined as in (\ref{F2dec}), we have
the vanishings: \begin{enumerate} \item $
H^i(G,\bigwedge^i\mathcal{F}^*)=0,\,\, n+2-i\geq0,\,i\geq1$ \item
$H^{i-1}(G,\bigwedge^i\mathcal{F}^*)=0,\,\,n+2-i=0 .$
\item $H^{i-2}(G,\bigwedge^i\mathcal{F}^*)=0,\,\,n+2-i=0 .$
\item $H^{i-3}(G,\bigwedge^i\mathcal{F}^*)=0,\,\,n+2-i\leq0,\,i\geq1  .$
\end{enumerate}
\end{lemm}
The proof of the Lemma is postponed to an appendix. \cqfd

As an immediate corollary, we get the following:
\begin{coro}\label{coro1}Under the same assumptions as in Proposition \ref{propvan}, the numbers
$$rk\,H^1(X,K_X\otimes\mathcal{I}_W(nH)),\,\,rk\,H^1(X,\mathcal{I}_W(nH))$$
are bounded by  $Cn^2$ for some constant $C$.
\end {coro}

 Combining  Corollary \ref{coro1}  with Riemann-Roch and Barth's Theorem
\ref{barth}, we get the following corollary:
\begin{coro} \label{c}The spaces $H^0(X,K_X(nH)\otimes \mathcal{I}_W)$ and $H^0(X,\mathcal{I}_W(nH))$
have dimension bounded by $cn^2$ for some constant $c$.
\end{coro}
We shall use the following consequence of the uniform position principle of Harris:
\begin{lemm} \label{coro2}Let $A$ be generic and let $W'\subset W$ be a subset of $W=Sing\,S_A$.
Then if
$$  H^0(X,K_X(nH)\otimes\mathcal{I}_{W'})\not=H^0(X,K_X(nH)\otimes\mathcal{I}_{W}),$$
 $W'$ imposes ${\rm card}\,W'$
independent conditions to $H^0(X,K_X(nH))$.
Similarly, if $$  H^0(X,\mathcal{O}_X(nH)\otimes\mathcal{I}_{W'})\not=
H^0(X,\mathcal{O}_X(nH)\otimes\mathcal{I}_{W}),$$
 $W'$ imposes ${\rm card}\, W'$
independent conditions to $H^0(X,\mathcal{O}_X(nH))$.
\end{lemm}
{\bf Proof.} Indeed, we represented in the previous proof the set $W$ as the projection in
$X$ of a $0$-dimensional subscheme $\widetilde{W}$ of
$G\times X$, defined as the zero set of a generic transverse section
of a very ample vector bundle $\mathcal{F}\boxtimes\mathcal{O}_X(H)$ on $G\times X$.
Thus the uniform position principle \cite{harris} applies to
$\widetilde{W}$, and allows to conclude that all subsets of $W$ of given cardinality
impose the same number of independent conditions to
$H^0(X,K_X(nH))$ or $H^0(X,\mathcal{O}_X(nH))$. This number is then obviously equal to
$$Min\,({\rm card}\, W',\, a)$$
where $a=rk\,(rest:\,H^0(X,K_X(nH)){\rightarrow} H^0(W,K_X(nH)_{\mid W})))$, resp.
$$a=rk\,(rest:\,H^0(X,\mathcal{O}_X(nH)){\rightarrow} H^0(W,\mathcal{O}_W(nH)))$$
in the second case.

\cqfd

From now on, we will treat separately the uniruled and the $K$-trivial cases.

\vspace{1cm}

{\bf The uniruled case.} We may assume $(X,H)$ satisfies the
inequality $H^2K_X<0$ of Lemma \ref{KH2}. We want to study the set
$Z$ of (\ref{Z}) for a generic surface $S\in \mid nH\mid$, and more
precisely the irreducible components $Z'$ of $Z$ which are of
dimension $\geq dim\,\mathbb{P}(V')$.

Degenerating
$S$ to $S_A$ and applying Lemma \ref{voisinlemme}, we find that
the specialization
$Z'_s$ of $Z'$ is contained in
$$
Z_0:=\{(v,v')\in \mathbb{P}(V_A)\times\mathbb{P}(V'_A),\,vv'_{\mid W}=0\},
$$
where $$V_A=H^0(S_A,K_{S_A}),\,V'_A=H^0(S_A,\mathcal{O}_{S_A}(nH)).$$
\begin{prop} \label{propspe}$Z'_s$ is contained in the union

\begin{eqnarray}
\label{union}
\mathbb{P}H^0(S_A,K_{S_A}\otimes\mathcal{I}_W)\times\mathbb{P}(V'_A)
\,\,\cup\,\,\mathbb{P}(V_A)\times\mathbb{P}H^0(S_A,\mathcal{O}_{S_A}(nH)\otimes\mathcal{I}_W)
.
\end{eqnarray}
\end{prop}
{\bf Proof.}
We observe that $Z_0$ is a union of
irreducible components indexed by subsets $W'\subset W$, with complementary set
$W'':=W\setminus W'$:
$$Z_0=\cup_{W'\subset W}Z_{W'},\,Z_{W'}:=\mathbb{P}H^0(S_A,K_{S_A}\otimes\mathcal{I}_{W'})\times
\mathbb{P}H^0(S_A,\mathcal{O}_{S_A}(nH)\otimes\mathcal{I}_{W''}).$$
We use now Lemma \ref{coro2}: it says  that if both conditions
$$H^0(X,K_X(nH)\otimes\mathcal{I}_{W'})\not=H^0(X,K_X(nH)\otimes\mathcal{I}_{W}),$$
$$ H^0(X,\mathcal{O}_X(nH)\otimes\mathcal{I}_{W''})\not=
H^0(X,\mathcal{O}_X(nH)\otimes\mathcal{I}_{W})$$ hold, then $W'$
imposes ${\rm card}\,W'$ independent conditions to $H^0(X,K_X(nH))$
and $W''$ imposes ${\rm card}\,W''$ independent conditions to
$H^0(X,\mathcal{O}_X(nH))$. Thus the codimension of $Z_{W'}$ in
$\mathbb{P}(V_A)\times\mathbb{P}(V'_A)$ is equal to ${\rm
card}\,W'+{\rm card}\,W''={\rm card}\,W$. But  ${\rm card}\,W$ is
equal to $\frac{n(n^2-1)}{6}H^3$ by Theorem \ref{barth},  while the
dimension of $V_A=H^0(S_A,K_{S_A})\cong H^0(X,K_X(nH))$ is equal to
$$\frac{1}{6}n^3H^3+\frac{1}{4}n^2K_XH^2+\,{\rm affine\,\,linear\,\,term\,\,in}\,\, n$$ by Riemann-Roch.

As $K_XH^2<0$, we conclude that for $n$ large enough, if $W'$ is as above,
we have $dim\,Z_{W'}<dim\,\mathbb{P}(V'_A)$.
Thus, for large $n$, the only components of $Z_0$ which may have dimension
$\geq dim\,\mathbb{P}(V'_A)$ are
$\mathbb{P}H^0(S_A,K_{S_A}\otimes\mathcal{I}_W)\times\mathbb{P}(V'_A)
$ and $\mathbb{P}(V_A)\times\mathbb{P}H^0(S_A,\mathcal{O}_{S_A}(nH)\otimes\mathcal{I}_W)$.
\cqfd
\begin{coro} \label{coroproj}
Assume $S$ is generic and $Z'\subset
\mathbb{P}(V)\times\mathbb{P}(V')$ is an irreducible component of
$Z$ which has dimension $\geq dim\,\mathbb{P}(V')$. Then either

i) $dim\,pr_1(Z')\leq cn^2$ or

ii) $dim\,pr_1(Z')\leq cn^2$,

where $c$ is the constant of Corollary \ref{c}.

\end{coro}

{\bf Proof.} By Proposition \ref{propspe}, the specialization $Z'_s$
of $Z'$ is contained in the union (\ref{union}). As we have by
Corollary \ref{c}
$$dim\,\mathbb{P}H^0(X,K_X(nH)\otimes \mathcal{I}_W)<cn^2,\,dim\,\mathbb{P}H^0(X,\mathcal{I}_W(nH))<cn^2,$$
 this implies that the cycle $Z'_s$ satisfies:
$$h_1^{nc^2}h_2^{nc^2}[Z'_s]=0\,\,{\rm in}\,\,H^*(\mathbb{P}(V_A)\times\mathbb{P}(V'_A),\mathbb{Z}),$$
where
$$h_1:=pr_1^*c_1(\mathcal{O}_{\mathbb{P}(V_A)}(1)),\,h_2:=pr_2^*c_1(\mathcal{O}_{\mathbb{P}(V'_A)}(1)),$$
and $[Z'_s]$ is the cohomology class of the cycle $Z'_s$.

It follows that we also have
\begin{eqnarray}h_1^{nc^2}h_2^{nc^2}[Z']=0\,\,{\rm
in}\,\,H^*(\mathbb{P}(V)\times\mathbb{P}(V'),\mathbb{Z}).
\label{poscond}
\end{eqnarray}
We claim that this implies that i) or ii) holds. Indeed, as $Z'$ is
irreducible of dimension $\geq 2cn^2$, there are well defined
generic ranks $k_1 ,\,k_2$ of the projection $pr_{1\mid
Z'},\,pr_{2\mid Z'}$ respectively, which are also the generic  ranks
of the pull-back of the $(1,1)$-forms $pr_1^*\omega_1$,
$pr_2^*\omega_2$ to $Z'$, where $\omega_i$ are the Fubini-Study
$(1,1)$-forms on $\mathbb{P}(V),\,\mathbb{P}(V')$. As the form
$$pr_1^*\omega_1^{nc^2}\wedge pr_1^*\omega_2^{nc^2}$$
is semi-positive on $Z'$, the condition (\ref{poscond}) implies that
everywhere on $Z$, we have
$$pr_1^*\omega_1^{nc^2}\wedge pr_2^*\omega_2^{nc^2}=0.$$
  As
$dim\,Z'\geq 2cn^2$ and  $(pr_1,pr_2)$ is an immersion on
the smooth locus of $Z'$, this implies easily that either $k_1=rk\,pr_1$ or
$k_2=rk\,pr_2$ has to be $<cn^2$, that is i) or ii). \cqfd
\begin{coro}\label{lastcoro} With the same assumptions as in the previous corollary, if $(v,v')\in Z'$,
one has either

i) $rk\,\mu_v:V'\rightarrow H^1(S,\Omega_S)< cn^2$, or

ii) $rk\,\mu_{v'}:V\rightarrow H^1(S,\Omega_S)< cn^2$,

where i) and ii) refer to the two cases of Corollary \ref{coroproj} and
$$\mu_v(\cdot)=\mu(v\otimes\cdot),\,\mu_{v'}(\cdot)=\mu(\cdot\otimes v').$$\end{coro}
{\bf Proof.} Indeed, assume case i) of Corollary \ref{coroproj} holds. As
$dim\,Z'\geq dim\,\mathbb{P}(V')$, the generic fibre of
$pr_1:Z'\rightarrow \mathbb{P}(V')$ has dimension $> dim\,\mathbb{P}(V')-cn^2$.
But the generic fibre is, by definition of $Z$, equal to
$\mathbb{P}(Ker\,\mu_v)$. Thus $rank\,\mu_v<cn^2$.

In case ii), we can do the same reasoning, as we have
$$dim\,Z'\geq dim\,\mathbb{P}(V')\geq dim\,\mathbb{P}(V).$$

\cqfd The proof that such a $Z'$ does not exist, and thus, the proof
of Proposition \ref{propmainsec2} in the uniruled case, concludes
now by the following two Propositions :
\begin{prop}\label{propcn2fini} Let $S\in\mid nH\mid$  be generic, with
$n$ large enough. Let $c$ be any positive constant. Then there
exists a constant $A$ such that the sets \begin{eqnarray}
\label{Gamma} \Gamma=\{v\in V,\,rk\,\mu_v<cn^2\},\\
\label{Gamma'} \Gamma'=\{v'\in V',\,rk\,\mu_{v'}<cn^2\},
\end{eqnarray}
both have dimension bounded by $A$.

\end{prop}
\begin{prop}\label{propfiniempty}
Let $A$ be any positive constant. Let $S\in\mid nH\mid$  be generic,
with $n$ large enough  (depending on $A$).  Then the set
$$ B=\{v\in V,\,rk\,\mu_v<A\}$$
reduces to $0$.
\end{prop}

Indeed, we know by Corollary \ref{lastcoro} that our set $Z'$ should
satisfy either $pr_1(Z)\subset \overline{\Gamma}$ (case i) ) or
$pr_2(Z)\subset \overline{\Gamma'}$ (case ii) ), where the
$\overline{\cdot}$ means projection to projective space. Thus by
Proposition \ref{propcn2fini}, one concludes that in case i),
$dim\,pr_1(Z')\leq A$ and in case ii), $dim\,pr_2(Z)\leq A$, where
$A$ does not depend on $n$.

In case ii), it follows that $dim\,Z\leq dim\,\mathbb{P}(V)+A$ and
as we have $dim\,\mathbb{P}(V)+A<dim\,\mathbb{P}(V')$ by Lemma
\ref{12declem}, this gives a contradiction.

In case i), it follows, arguing  as in the proof of Corollary
\ref{lastcoro}, that for $(v,v')\in Z'$, one has $rk\,\mu_v<A$. This
is impossible unless $Z'$ is empty by Proposition
\ref{propfiniempty}. Thus, assuming Propositions \ref{propcn2fini}
and \ref{propfiniempty}, Proposition \ref{propmainsec2} is proved
for uniruled threefolds with $H^2(X,\mathcal{O}_X)=0$. \cqfd
\vspace{0,5cm}

{\bf Proof of Proposition \ref{propcn2fini}.} Our first step is to
reduce the statement to the case where $S$ is a surface in
$\mathbb{P}^3$. This is done as follows: we choose once for all a
morphism $$f:X\rightarrow\mathbb{P}^3$$ given by $4$ sections of
$H$, so that $f^* \mathcal{O}_{\mathbb{P}^3}(1)=H$. We shall prove
the result for surfaces of the form $S=f^{-1}(\Sigma)$, where
$\Sigma$ is a generic smooth surface of degree $n$ in
$\mathbb{P}^3$. Let $f_S:S\rightarrow\Sigma$ be the restriction of
$f$ to $S$. We have  trace maps
$$f_{S*}:H^1(S,\Omega_S(sH))\rightarrow
H^1(\Sigma,\Omega_\Sigma(s)),$$
$$ f_{S*}:H^0(S,K_S(sH))\rightarrow H^0(\Sigma,K_\Sigma(s))$$
for all integers $s$. We note now that the map $\mu$ admits obvious
twists that we shall also denote by $\mu$:
$$\mu:H^0(S,K_S(lH))\otimes H^0(S,\mathcal{O}_S(nH)) \rightarrow
H^1(S,\Omega_S(lH)).$$ Furthermore, we have  similarly defined
bilinear maps $\mu^\Sigma$:
$$\mu^\Sigma:H^0(\Sigma,K_\Sigma(l))\otimes H^0(\Sigma,\mathcal{O}_\Sigma(n)) \rightarrow
H^1(\Sigma,\Omega_\Sigma(l)).$$

All the maps $\mu$ can be defined using the maps
$$\delta:H^0(S,K_S(lH))\hookrightarrow H^1(S,\Omega_S((-n+l)H)),$$
induced by the exact sequence (which is itself a twist of the normal
exact sequence)
$$0\rightarrow \Omega_S(-nH)
\rightarrow {\Omega_X^2}_{\mid S}\rightarrow K_S\rightarrow0,$$
twisted by $lH$, and then the product map $$
H^1(S,\Omega_S((-n+l)H))\otimes H^0(S,\mathcal{O}_S(nH))\rightarrow
H^1(S,\Omega_S(lH)).$$ The same is true for the maps $\mu_\Sigma$.

 As
there is a commutative diagram of normal exact sequences
$$\begin{matrix}&0&\rightarrow &T_S&
\rightarrow& {T_X}_{\mid S}&\rightarrow&
\mathcal{O}_S(nH)&\rightarrow&0\\
&&&f_{S*}\downarrow&&f_*\downarrow&& \parallel&&\\
&0&\rightarrow &f^*T_\Sigma &\rightarrow
&f^*{T_{\mathbb{P}^3}}_{\mid S}&\rightarrow&
\mathcal{O}_S(nH)&\rightarrow&0
\end{matrix},
$$
where the bottom line is the normal bundle sequence of $\Sigma$
pulled-back to $S$, it follows  that for $v\in H^0(S,K_S(l))$ and
$\eta\in H^0(\Sigma,\mathcal{O}_\Sigma(n))$, we have:

\begin{eqnarray}
\label{10dec}f_{S*}(\mu_v(f_S^*\eta))=\mu^\Sigma_{f_{S*}(v)}(\eta),
\end{eqnarray}
 Equation (\ref{10dec}) implies that
$$rk\,(\mu^\Sigma_{f_{S*}(v)}:H^0(\Sigma,\mathcal{O}_\Sigma(n))\rightarrow
H^1(\Sigma,\Omega_\Sigma(l)))$$
$$\leq
rk\,(\mu_v:H^0(S,\mathcal{O}_S(n))\rightarrow H^1(S,\Omega_S(l))).
$$
Let us now prove the first case of Proposition \ref{propcn2fini},
namely for the set $\Gamma$. The second proof is done similarly.

Starting from a sufficiently ample $H$, one finds that
$H^0(X,\mathcal{O}_X(4H))$ restricts surjectively onto $
H^0(X_u,\mathcal{O}_{X_u}(4H))$, for any $u\in  \mathbb{P}^3$, where
$X_u:=f^{-1}(u)$.

We have the following Lemma:
\begin{lemm}\label{postponed} The image $\Gamma_\Sigma$ of the composed map
$$ \Gamma \times H^0(X,\mathcal{O}_X(4H))\stackrel{\nu}{\rightarrow} H^0(S,K_S(4H))
\stackrel{f_{S*}}{\rightarrow} H^0(\Sigma,K_\Sigma(4)),$$ where
$\nu$ is the product,  has dimension at least equal to
$\frac{1}{N}dim\,\Gamma$, where
$$N:=rk\,H^0(X,\mathcal{O}_X(4H)).$$
 \end{lemm}
{\bf Proof.} Indeed, as the restriction map
$H^0(X,\mathcal{O}_X(4H))\rightarrow H^0(X_u,\mathcal{O}_{X_u}(4H))$
is surjective,   if $e_i$ is a basis of $H^0(X,\mathcal{O}_X(4H))$,
the map
$$\Gamma\rightarrow \Gamma_\Sigma^D,$$
$$\gamma\mapsto f_{S*}(\gamma e_i),$$
is injective. Thus $dim\,\Gamma\leq Ndim\,\Gamma_\Sigma$.

\cqfd

 On the other hand, if $v\in
\Gamma,\,\alpha\in H^0(X,\mathcal{O}_X(4H))$, we have
$$rk\,\mu_{\alpha v}\leq rk\,\mu_{ v}$$
because $\mu_{\alpha v}=\alpha\mu_v$. Thus we conclude that the
following hold:
 $$dim\,\Gamma_\Sigma\geq \frac{1}{N}dim\,\Gamma,$$
$$rk\,\mu^\Sigma_w\leq rk\,\mu_v\leq cn^2,$$
for all $w\in \Gamma_\Sigma$.

As $N$ does not depend on $n$,  it suffices to show the result for
generic $\Sigma$ in $\mathbb{P}^3$ and for the product
$$\mu^\Sigma:H^0(\Sigma,K_\Sigma(4))\times
H^0(\Sigma,\mathcal{O}_\Sigma(n)) \rightarrow
H^1(\Sigma,\Omega_\Sigma(4)).$$ This last product is well known (cf
\cite{voisinbook},6.1.3) to identify to the multiplication in the
Jacobian ring of $\Sigma$:
$$\mu_\Sigma:H^0(\Sigma,\mathcal{O}_\Sigma(n))\times H^0(\Sigma,\mathcal{O}_\Sigma(n))
\rightarrow R^{2n}_\Sigma.$$ Thus we have to show that for generic
$\Sigma$, the set
$$\Gamma_\Sigma:=\{v\in H^0(\Sigma,\mathcal{O}_\Sigma(n)),\,rk\,\mu_{\Sigma,v}\leq cn^2\}
$$ has dimension bounded by a constant which is  independent of
$n$.

 For this, we specialize to the case where
$\Sigma$ is the Fermat surface, that is, its defining equation is
$\sigma=\sum_{0}^{3}X_i^n$. The Jacobian ideal of $\Sigma$ is then
generated by the $X_i^{n-1}$, and there is thus a natural action of
the torus $(\mathbb{C}^*)^4$ on the Jacobian ring $R_\Sigma$, by
multiplication of the coordinates by a scalar. The subspace
$$\Gamma_\Sigma\subset R^{n-1}_\Sigma$$ is thus invariant under
$(\mathbb{C}^*)^4$. Note that the fixed points of the induced action
on $\mathbb{P}(R^{n-1}_\Sigma)$ are the monomials, and are thus
isolated. It follows that we have the inequality
$$dim\,\overline{\Gamma}_\Sigma\leq \,\,{\rm number\,\,\,of\,\,fixed\,\,points\,\,on}\,\,\overline{\Gamma}_\Sigma.$$
Thus we have to bound the number of monomials
$$X_I=X_0^{i_0}X_1^{i_1}X_2^{i_2}X_3^{i_3},\,i_0+i_2+i_3+i_4=n,$$
such that
$$rk\,X_I:R^{n}_\Sigma\rightarrow R^{2n-1}_\Sigma\leq cn^2.$$
But the kernel of the multiplication by $X_I$ above is equal to the
ideal
$$ X_0^{n-i_0}S^{i_0}+\ldots X_3^{n-i_3}S^{i_3},$$
where $S^l:=H^0(\mathbb{P}^3,\mathcal{O}_{\mathbb{P}^3}(l))$,
and thus  has  dimension $\leq \sum_{k}rk\,S^{i_k}$. Hence, if
$rk\, X_I\leq cn^2$, we must have
\begin{eqnarray}\label{rang}\sum_{k}rk\,S^{i_k}\geq rk\,S^n-cn^2,
\end{eqnarray}
with $\sum_k i_k=n$. It is not hard to see  that there exists an
integer $l>0$ such that, if $n$ is large enough and (\ref{rang})
holds for $I,\,n$, one of the $i_k's$ has to be $\geq n-l$. Thus the
other $i_j$'s have to be non greater than $l$. This shows
immediately that the number of such monomials is bounded by a
constant independent of $n$ and concludes the proof.

 \cqfd
 {\bf Proof of Proposition \ref{propfiniempty}.} The key point
is the following fact from \cite{green}. \begin{prop} \label{green}
Let $X$ be any projective manifold and $H$ be a very ample  line
bundle on $X$. Let $A$ be a given constant, and for $n>A$,
let $M\subset H^0(X,\mathcal{O}_X(nH))$ be a subspace of
codimension $\leq A$. Then
$$H^0(X,\mathcal{O}_X(H))\cdot M\subset H^0(X,\mathcal{O}_X((n+1)H))$$ has codimension $\leq A$,
with strict inequality if $M$ has no base-point.
\end{prop}
Assume $v\in V$ satisfies the condition that $rk\,\mu_v<A$. Let
$M:=Ker\,\mu_v\subset H^0(S,\mathcal{O}_S(nH))$. By
Proposition \ref{green}, we conclude that if $n>A$, we have
$$H^0(S,\mathcal{O}_S(H))\cdot M\subset
H^0(S,\mathcal{O}_S((n+1)H))$$
has codimension $< A$. Next, we consider for each $l$ the map
$${\mu}_v^l:H^0(S,\mathcal{O}_S((n+l)H))\rightarrow
H^1(S,\Omega_S(l)),$$ obtained as the composite of the twisted Kodaira-Spencer map
$$H^0(S,\mathcal{O}_S((n+l)H))\rightarrow H^1(S,T_S(l)),$$
and the contraction with $v$, using the contraction map
$$ H^0(S,K_S)\otimes H^1(S,T_S(l))\rightarrow H^1(S,\Omega_S(l)).$$
We note that the kernel $M_l$ of the map ${\mu}_v^l$ contains
$$M_1\cdot H^0(S,\mathcal{O}_S((l-1)H)).$$
On the other hand, $M_1$ also contains the image of the map
$$ H^0(S,T_X(H)_{\mid S})\rightarrow H^0(S,\mathcal{O}_S((n+1)H))$$
induced by the normal bundle sequence twisted by $H$. We may assume
that $H$ is ample enough so that $H^0(X,T_X(1))$ is generated by
global sections, and then $M_1$ has no base-point. Proposition
\ref{green}  thus implies that if $n>A$, the numbers
$$corank\,M_l$$
are strictly decreasing, starting from $l\geq1$.  Hence we conclude
that $$M_A=H^0(S,\mathcal{O}_S((n+A)H)).$$ As $n$ is large and $A$
is fixed, we may assume that
$$ H^0(X,K_X((2n-A)H))\otimes H^0(X,\mathcal{O}_X((n+A)H))
\rightarrow H^0(X,K_X(3nH))$$ is surjective, and that the same is true after
restriction to $S$. Thus we conclude that
$$M_A\cdot H^0(S,K_X((2n-A)H)_{\mid S})=H^0(S,K_X(3nH)_{\mid S}).$$
 We use now the definition of $M_{A}$, and the compatibility of the
 twisted Kodaira-Spencer  maps and the maps
$\lrcorner v$ with multiplication. This implies that for any $P\in
H^0(S,K_X((3n)H)_{\mid S})$, sending to $\overline{P}\in
H^1(S,T_S(K_S(2nH)))$, via  the
map induced by the twisted normal bundle sequence $$ 0\rightarrow
T_S(K_S(2nH))\rightarrow {T_X}_{\mid S}(K_S(2nH))\rightarrow
K_X(3nH)_{\mid S}\rightarrow0,$$ we have
\begin{eqnarray}\label{premiervan}
\overline{P}\lrcorner v=0\,\,{\rm in}\,\, H^1(S,\Omega_S(K_S(nH))).
\end{eqnarray}
 We have now a map
 $$\delta:H^1(S,\Omega_S(K_S(nH)))\rightarrow H^2(S,K_S),$$
 induced by the exact sequence
 $$0\rightarrow K_S\rightarrow \Omega_X(K_X(2nH))_{\mid
S}\rightarrow \Omega_S(K_S(nH))\rightarrow0,$$ and one knows (cf
\cite{carlsongriffiths}) that up to a multiplicative coefficient,
one has
\begin{eqnarray}
\label{eq5dec}\delta(\overline{P}\lrcorner v) =<v, res_S(P)>,
\end{eqnarray} where on the right, $<,>$ is Serre duality between
$H^0(S,K_S)$ and $H^2(S,\mathcal{O}_S)$, and the Griffiths  residue
map \begin{eqnarray}
\label{residue}H^0(X,K_X(3nH))\stackrel{res_S}{\rightarrow}
H^2(S,\mathcal{O}_S) \end{eqnarray} is described in
\cite{voisinbook},6.1.2.
  The key point for us is
that, because in our case $H^3(X,\Omega_X)=0$ and  because $n$ is
large enough,  the residue map (\ref{residue}) is surjective, and
thus (\ref{premiervan}) together with (\ref{eq5dec}) imply that, for
all $\eta\in H^2(S,\mathcal{O}_S)$, one has
$$<\eta,v>=0,$$
which implies that $v=0$.
 \cqfd {\bf The
Calabi-Yau case.} Here $X$ has trivial canonical bundle and
satisfies $H^2(X,\mathcal{O}_X)=0$. We use in this case a variant of
Lemma \ref{lemdim}. As $K_X$ is trivial, the spaces $V$ and $V'$ are
equal, and the pairing $\mu:V\times V'\rightarrow H^1(S,\Omega_S)$
is symmetric. Thus, using Bertini,  Lemma \ref{lemdim} can be
refined as follows (cf \cite{voisinintjmath}):
\begin{lemm}\label{lemdimsym} Let $\mu:V\otimes V'\rightarrow H^1(S,\Omega_S)$
be symmetric and
$q:H^1(S,\Omega_S)\rightarrow S^2V^*$ be its dual.
Then the generic quadric
in $Im\,q$ is non-singular if the following condition holds. There is no
subset
$Z\subset\mathbb{P}(V)$ contained in
the base-locus of $Im\,q$ and satisfying:
$$rk\,\mu_v\leq dim\,Z,\,\forall v\in  Z.$$

\end{lemm}

We have to verify that such a $Z$ does not exist for generic $S\in \mid nH\mid$, $n$ large enough.
Degenerating $S$ to $S_A$ as before, the base-locus of $Im\,q$ specializes to
a subspace of the base-locus of
$Im\,q_A$. We now use Lemma \ref{voisinlemme}, together with
Corollary \ref{c}, to conclude that the base-locus of
$Im\,q_A$ has dimension $\leq cn^2$, for some $c$ independent of $n$.

Thus the base-locus of $Im\,q$ also has dimension bounded by
$cn^2$, for generic $S$.

By definition of $Z$, it follows that for
$v\in Z$ one has
$$rk\,(\mu_v:V\rightarrow H^1(S,\Omega_S))\leq cn^2.$$

Using Proposition \ref{propcn2fini}, it follows that $dim\,Z\leq A$ for some constant $A$ independent of
$N$. But then,
for
$v\in Z$, one has
$$rk\,(\mu_v:V\rightarrow H^1(S,\Omega_S))\leq A$$
which implies that $Z$ is empty by Proposition \ref{propfiniempty}.
This concludes the proof of Proposition \ref{propmainsec2} when $X$
is a Calabi-Yau threefold. \cqfd

\section{The case where $H^2(X,\mathcal{O}_X)\not=0$ \label{sec3}}
In this section, we show how to adapt
the previous proof to the case where $X$ is uniruled with $H^2(X,{\mathcal
O}_X)\not=0$.

In this case, a smooth  birational model of $X$
admits a map $\phi:X'\dasharrow \Sigma$, with generic
fibre isomorphic to ${\mathbb P}^1$, where $\Sigma$ is a smooth surface.
Note that  $\phi_*$ sends $H^3(X,\Omega_X)$ isomorphically
to $H^2(\Sigma,{\mathcal O}_\Sigma)$.

We may assume that $X'$ carries a line bundle $H$ such that
$$H^2K_{X'}<0,$$
because there is a smooth birational model of $X$ on which such an
$H$  exists, and by blowing-up this $X'$ to an $\widetilde{X}$ with exceptional
relatively anti-ample divisor
$E$,  we
may assume that $\phi$ becomes defined, while an $\widetilde{H}$ of
the form $\tau^*H-\epsilon E$ with small $\epsilon$ will still
satisfy   the property $\widetilde{H}^2.K_{\widetilde{X}}<0$.

In the sequel $X,\,H,\,\phi$ will satisfy the properties above.
For $S$ a smooth surface in $\mid nH\mid$,
we have the Gysin maps:
$$\phi_*:H^1(S,\Omega_S)\rightarrow H^1(\Sigma,\Omega_\Sigma),
\,\phi_*:H^2(S,{\mathcal O}_S)\rightarrow H^2(\Sigma,{\mathcal
O}_\Sigma), $$
$$\phi_*:H^2(S,{\mathbb Z})\rightarrow
H^2(\Sigma,{\mathbb Z}).$$ We will denote by
$$H^1(S,\Omega_S)_\Sigma,\,H^2(S,{\mathcal O}_S)_\Sigma,\,
H^2(S,{\mathbb Z})_\Sigma$$
the respective kernels of these maps.
The proof will use the following variant of Proposition
\ref{prop1}:
\begin{prop}\label{propvar} Assume there is a $S\in \mid nH\mid$, and a
$\lambda\in H^1(S,\Omega_S)_\Sigma$ such that the natural map
$$\mu_\lambda :H^0(S,{\mathcal O}_S(nH))
\rightarrow H^2(S,{\mathcal O}_S)_\Sigma$$
defined as in (\ref{mulambda}) is surjective.
Then the Hodge conjecture is true for integral Hodge classes
on $X$.
\end{prop}
{\bf Proof.}
We consider a simply connected
open set in $\mid nH\mid$ parameterizing
smooth surfaces and containing the point $0\in \mid nH\mid$
parametrizing $S$.
We study the infinitesimal variation of Hodge structure on
$H^2(S_t,{\mathbb Z})_\Sigma$ for $t\in B$.

By the same reasoning as in the proof of Proposition \ref{prop1},
the existence of $\lambda$ satisfying the property above implies
that at some point $\lambda\in H^{1,1}(S)_{{\mathbb R},\Sigma}$, the
natural map

$$\psi:{ H}^{1,1}_{{\mathbb R},\Sigma}
\rightarrow H^2(S,{\mathbb R})_\Sigma$$
is a
submersion. Here on the left hand
side, we have the real vector bundle with fibre
$H^{1,1}(S_t)_{{\mathbb R},\Sigma}$ at the point $t$, and  on each
fibre $H^{1,1}(S_t)_{{\mathbb R},\Sigma},\,\psi$ is the inclusion
$H^{1,1}(S_t)_{{\mathbb R},\Sigma} \subset H^2(S_t,{\mathbb
R})_\Sigma$, followed by the topological isomorphism
$H^2(S_t,{\mathbb R})_\Sigma\cong H^2(S,{\mathbb R})_\Sigma$.

This implies that the image of $\psi$ contains an open cone and we
deduce from this as in the proof of Proposition \ref{prop1} that
$H^2(S,{\mathbb Z})_\Sigma$ is generated over ${\mathbb Z}$ by
classes $\alpha$ which are algebraic on some nearby fiber $S_t$.

Consider now the inclusion $j:S\rightarrow X$. It induces a surjective
Gysin map
$$H^2(S,{\mathbb Z})\rightarrow H^4(X,{\mathbb Z})$$
by Lefschetz hyperplane theorem.
On the other hand, we have a commutative diagram of
Gysin maps:
$$\begin{matrix}&H^2(S,{\mathbb Z})&\stackrel{j_*}{\rightarrow}& H^4(X,{\mathbb
Z})&\cr
&\phi_*\downarrow&&\phi_*\downarrow&\cr
&H^2(\Sigma,{\mathbb Z})&=&H^2(\Sigma,{\mathbb Z})&.
\end{matrix}
$$
From this and the previous conclusion,  we deduce  that the group
$$Ker\,(\phi_*:H^4(X,{\mathbb Z})\rightarrow H^2(\Sigma,{\mathbb Z}))
=j_*H^2(S,\mathbb{Z})_\Sigma
$$
is generated by classes of algebraic cycles on $X$.

Thus it suffices to prove the following:
\begin{lemm}\label{lemmeSigma}
Let $\alpha$ be an integral Hodge class of degree $2$ on $\Sigma$
which is in  $$Im\,(\phi_*:H^4(X,{\mathbb Z})
\rightarrow H^2(\Sigma,{\mathbb Z})).$$
Then there is an algebraic $1$- cycle $Z$ on $X$ such that
$$\alpha=\phi_*([Z]).$$
\end{lemm}
Indeed, assuming the Lemma, if $\alpha$ is an integral  Hodge class
on $X$ of degree $4$, $\phi_*\alpha$ is an integral Hodge class of
degree $2$ on $\Sigma$, hence  is equal to $\phi_*([Z])$ for some $Z$.
Hence $\alpha-[Z]$ belongs to $Ker\,\phi_*$ and thus it is algebraic
as we already proved. This proves the Proposition.\cqfd {\bf Proof
of Lemma \ref{lemmeSigma}.} We may assume by Lefschetz $(1,1)$
theorem and because $\Sigma$ is algebraic, that $\alpha$ is the
class of a curve $C\subset S$ which is in general position. Thus
$$\phi_C: X_C:=\phi^{-1}(C)\rightarrow C$$
is a geometricall ruled surface, which admits a section $C'\subset X_C$
(see \cite{beauville}, or \cite{harrismazur} for a more general statement).

But then the curve $C'\subset X$ satisfies
$\phi_*[C']=[C]$.
\cqfd

By Proposition \ref{propvar}, the proof of Theorem \ref{main} in
case where $X$ is uniruled and satisfies
$H^2(X,\mathcal{O}_X)\not=0$ will now be  a consequence of the
following proposition.
\begin{prop} \label{propmainsec3}Let the pair $(X,H)$ satisfy the inequality
$H^2K_X<0$. Then for $n$ large enough, for $S$ a generic surface in
$\mid nH\mid$, there is a $\lambda\in H^1(S,\Omega_S)_\Sigma$ which
satisfies the property that
$$\mu_\lambda:H^0(S,\mathcal{O}_S(nH))\rightarrow
H^2(S,\mathcal{O}_S)_\Sigma$$ is surjective.
\end{prop}

The proof works exactly as the proof of Proposition
\ref{propmainsec2} in the uniruled case. The only thing to note is the fact that the analogue of
Proposition \ref{propfiniempty} still holds in this case, with
$V=H^0(S,K_S)_\Sigma$, $V'=H^0(S,\mathcal{O}_S(nH))$. This is indeed the only place where we used
the assumption $H^2(X,\mathcal{O}_X)=0$.

In this case, we have an isomorphism
$$\phi_*:H^3(\Omega_X)\cong H^2(\Sigma,\mathcal{O}_\Sigma),$$
so that for $S\subset X$ a smooth surface
$$H^2(S,\mathcal{O}_S)_\Sigma=Ker\,(j_*:H^2(S,\mathcal{O}_S)\rightarrow
H^3(X,\Omega_X)),$$ where $j$ is the inclusion of $S$ into $X$.

But the theory of Griffiths residues shows that the last kernel is
precisely generated by residues $res_S\omega,\,\omega\in
H^0(X,K_X(3nH))$. Thus, the arguments of Proposition
\ref{propfiniempty} will show in this case that if $v\in H^0(S,K_S)$
satisfies $rank\,\mu_v\leq A$, where $A$ is a given constant, and
$S\in \mid nH\mid$ with $n$ large enough, then
$$v\in (Ker\,j_*)^\perp,$$
where $\perp$ refers to Serre duality between $H^0(S,K_S)$ and
$H^2(S,\mathcal{O}_S)$. But as $Ker\,\phi_*=Ker\,j_*$, we have
$$(Ker\,j_*)^\perp=\phi^* H^0(\Sigma,K_\Sigma).$$
Thus if furthermore $v\in H^0(S,K_S)_\Sigma$, we must have $v=0$
because $$H^0(S,K_S)_\Sigma\cap \phi^* H^0(\Sigma,K_\Sigma)=0.$$ \cqfd

\section{Appendix}
We give for the convenience of the reader the proof of the vanishing
Lemma \ref{vanlemm3dec}. Recall that we want to prove the vanishing
of the spaces:
\begin{enumerate} \item\label{i} $
H^i(G,\bigwedge^i\mathcal{F}^*),\,\, n+2-i\geq0,\,i\geq1$
\item\label{ii} $H^{i-1}(G,\bigwedge^i\mathcal{F}^*),\,\,n+2-i=0 .$
\item\label{iii} $H^{i-2}(G,\bigwedge^i\mathcal{F}^*),\,\,n+2-i=0 .$
\item\label{iv} $H^{i-3}(G,\bigwedge^i\mathcal{F}^*),\,\,n+2-i\leq0  .$
\end{enumerate}

We use first the dual of the exact sequence (\ref{exactF}) to get a
resolution of $\bigwedge^i\mathcal{F}^*$:
$$\ldots\rightarrow   \bigwedge^{i-1}
(K\otimes\mathcal{E}^*) \otimes\mathcal{L}^{-1} \rightarrow
\bigwedge^i
(K\otimes\mathcal{E})^*\rightarrow\bigwedge^i\mathcal{F}^*
\rightarrow0.$$ This induces a spectral sequence converging to
$$H^i(G,\bigwedge^i\mathcal{F}^*),\, H^{i-1}(G,\bigwedge^i\mathcal{F}^*)
,\,H^{i-2}(G,\bigwedge^i\mathcal{F}^*),\,H^{i-3}(G,\bigwedge^i\mathcal{F}^*),$$
whose $E_1$ terms are
$$ {\rm Case \,\ref{i} }\,\,\,\,\,\,\,\,\,\,\,\,\,\,\,\,H^{i+s}(G,\bigwedge^{i-s}
(K\otimes\mathcal{E}^*) \otimes\mathcal{L}^{-s}),\,n+2\geq
i\geq1,\,i\geq s\geq0,$$
$${\rm Case  \,\ref{ii} }\,\,\,\,\,\,\,\,\,\,\,\,\,\,\,\,H^{i+s-1}(G,\bigwedge^{i-s}
(K\otimes\mathcal{E}^*) \otimes\mathcal{L}^{-s}),\, i=n+2,\,i\geq
s\geq0,$$
$${\rm Case  \,\ref{iii} }\,\,\,\,\,\,\,\,\,\,\,\,\,\,\,\, H^{i+s-2}(G,\bigwedge^{i-s}
(K\otimes\mathcal{E}^*) \otimes\mathcal{L}^{-s}),\, i=n+2\,i\geq
s\geq0$$
$${\rm Case \, \ref{iv} }\,\,\,\,\,\,\,\,\,\,\,\,\,\,\,\, H^{i+s-3}(G,\bigwedge^{i-s}
(K\otimes\mathcal{E}^*) \otimes\mathcal{L}^{-s}),\, n+2\leq
i,\,i\geq s\geq0$$ respectively.

Let $P\subset \mathbb{P}(K)\times G$ be the incidence scheme, so $P$
is a $\mathbb{P}^1$-bundle over $G$. Let $pr_i,\,i=1,\,2$ denote the projections from
$P$ to $\mathbb{P}(K)$ and $ G$ respectively. Let $H:=pr_1^*\mathcal{O}(1)$
and denote  also by $\mathcal{L}$ the pull-back of $\mathcal{L}$ to
$P$. Then $pr_2^*\mathcal{E}^*$ fits into an exact sequence:
$$0\rightarrow H^{-1}\rightarrow pr_2^*\mathcal{E}^* \rightarrow H\otimes \mathcal{L}^{-1}
\rightarrow0.$$ Thus the bundle
$$pr_2^*(\bigwedge^{i-s}
(K\otimes\mathcal{E}^*) \otimes\mathcal{L}^{-s})$$ admits a
filtration  whose successive quotients  are line bundles of the form
$$H^{-\alpha}\otimes(H\otimes \mathcal{L}^{-1})^{\beta}\otimes\mathcal{L}^{-s}=
H^{-\alpha+\beta}\otimes\mathcal{L}^{-\beta-s},$$ where
$\alpha+\beta=i-s,\,\alpha\geq0,\,\beta\geq0$. As we are interested
in
$$H^*(G,\bigwedge^{i-s}
(K\otimes\mathcal{E}^*) \otimes\mathcal{L}^{-s})=
H^*(G,R^0pr_{2*}(pr_2^*(\bigwedge^{i-s} (K\otimes\mathcal{E}^*)
\otimes\mathcal{L}^{-s}))),$$ it suffices to study the cohomology
groups
$$H^*(P,H^{-\alpha+\beta}\otimes\mathcal{L}^{-\beta-s}),$$
with $-\alpha+\beta\geq0$. These groups are equal to the groups
$$ H^*(G,S^{-\alpha+\beta}\mathcal{E}\otimes\mathcal{L}^{-\beta-s})$$
which are partially computed in \cite{voisinjems}. The conclusion is
the following:
\begin{lemm}\label{lemjems}

a) These groups  vanish for $*\not=n-2,\,2(n-2)$ and for
$\beta+s\leq n-2$.

b) For $*=n-2$, these groups  vanish if $-s-\alpha+1<0$.

c) For $*=2(n-2)$, these groups vanish if $-s-\alpha\geq -n+1 $.

\end{lemm}

{\bf Case \ref{i}.}  Here $*=i+s$, and the following inequalities
hold: \begin{eqnarray} \label{ineq3dec} \beta\geq\alpha\geq0,\,
\beta+s\geq n-1 \end{eqnarray} and furthermore
$$1\leq i\leq n+2,\,\alpha+\beta=i-s.$$
According to Lemma \ref{lemjems}, in order to get a non trivial
cohomology group, we have only two possibilities:

a) $i+s=n-2,\,-s-\alpha+1\geq0.$

b) $i+s=2(n-2),\,-s-\alpha<-n+1$.

In case a), we have $\beta+s\geq n-1$ and $\alpha+\beta+2s=i+s=n-2$,
which is clearly a contradiction as $\alpha+s\geq0$.

In case b), we have $\beta+s\geq n-1,\,\alpha+s\geq n$ and thus
$$2n-1\leq\alpha+\beta+2s=i+s=2(n-2)$$
which is clearly a contradiction.
\cqfd

{\bf Case \ref{ii}.}
 Now $*=i+s-1$ and $i=n+2$. We have again the inequalities
 (\ref{ineq3dec}) and furthermore
 $$i= n+2,\,\alpha+\beta=i-s.$$
By Lemma \ref{lemjems}, in order to get a non trivial cohomology
group, we have only two possibilities:

a) $i+s-1=n-2,\,-s-\alpha+1\geq0.$

b) $i+s-1=2(n-2),\,s+\alpha\geq n$.

In case a), we have $i=n+2$ and $s\geq0$, hence $i+s-1=n-2$ is
impossible.

In case b), we have $i+s=2n-3$, while $s+\alpha\geq n$ and
$\beta+s\geq n-1$ give $\alpha+\beta+2s=i+s\geq2n-1$, contradiction.

 \cqfd {\bf Case \ref{iii}.}
Now $*=i+s-2$ and $i=n+2$. We have again the inequalities
 (\ref{ineq3dec}) and furthermore
 $i= n+2,\,\alpha+\beta=i-s$.
As before, in order to get a non trivial cohomology group, we have
only two possibilities:

a) $i+s-2=n-2,\,-s-\alpha+1\geq0.$

b) $i+s-2=2(n-2),\,s+\alpha\geq n$.

In case a), we have $i=n+2$ and $s\geq0$, hence $i+s-2=n-2$ is
impossible.

In case b), we have $i+s=2n-2$, while $s+\alpha\geq n$ and
$\beta+s\geq n-1$ give $\alpha+\beta+2s=i+s\geq2n-1$, contradiction.

 \cqfd

{\bf Case \ref{iv}.} Now $*=i+s-3$ and $i\geq n+2$. We have again the
inequalities
 (\ref{ineq3dec}) and furthermore
 $i\geq n+2,\,\alpha+\beta=i-s$.
As before, in order to get a non trivial cohomology group, we have
only two possibilities:

a) $i+s-3=n-2,\,-s-\alpha+1\geq0.$

b) $i+s-3=2(n-2),\,s+\alpha\geq n$.

In case a), we have $i\geq n+2$ and $s\geq0$ thus $i+s-3=n-2$  is
impossible.

In case b), we have $i+s=2n-1$, while $s+\alpha\geq n$ and
$\beta+s\geq n-1$ give $\alpha+\beta+2s=i+s\geq2n-1$. Thus we must
have the two equalities $$ s+\alpha= n,\,\,\beta+s= n-1.$$ This
contradicts the fact that $\beta\geq\alpha$.

\cqfd

\end{document}